\documentclass[11pt,fullpage]{article}

\usepackage[pdftex]{graphicx}
\usepackage{latexsym}
\usepackage{color}
\usepackage{multirow}

\usepackage{amsmath}
\usepackage{amssymb}
\usepackage{amsthm}
\usepackage{amsfonts}
\usepackage{bm}

\newcommand{\fr}[1]{(\ref{#1})}  

\renewcommand{\baselinestretch} {1.3}

\makeatletter 
\def\singlespace{\def\baselinestretch{1}\@normalsize}

\newtheorem{theorem}{Theorem}
\newtheorem{lemma}{Lemma}

\newtheorem{proposition}{Proposition}

\newcommand{\bfm}[1]{\mbox{\boldmath $#1$}}
\newcommand{\bbeta}{\bfm{\beta}}

\newcommand{\bx}{{\bf a}}

\newcommand{\cA}{{\cal A}}

\newcommand{\Var}{\mbox{Var}}
\newcommand{\Cov}{\mbox{Cov}}

\long\def\ignore#1{}

\newcommand{\be}{\begin{equation}}
\newcommand{\ee}{\end{equation}}
\newcommand{\beqn}{\begin{eqnarray}}
\newcommand{\eeqn}{\end{eqnarray}}
\newcommand{\bes}{\begin{equation*}}
\newcommand{\ees}{\end{equation*}}
\newcommand{\beqns}{\begin{eqnarray*}}
\newcommand{\eeqns}{\end{eqnarray*}}

\newcommand{\lkr}{\left(} 
\newcommand{\lkv}{\left[} 
\newcommand{\rkv}{\right]} 
\newcommand{\rkr}{\right)}  
\newcommand{\lfi}{\left\{}  

\begin{document}


 \title{\bf Classification with many classes: challenges and pluses }


\author{{\bf Felix Abramovich}\\
Department of Statistics\\
 and Operations Research \\
Tel Aviv University \\
Tel Aviv 69978 \\
Israel \\
{\it felix@post.tau.ac.il}
\and
{\bf Marianna Pensky} \\
Department of Mathematics\\
University of Central Florida \\
4393 Andromeda Loop N \\ 
Orlando, FL 32816 \\
USA \\
{\it Marianna.Pensky@ucf.edu}
}

\date{}

\maketitle

\begin{abstract}
The objective of the paper is to study accuracy of multi-class classification
in high-dimensional setting, where the number of classes is also large (``large $L$, large $p$, small $n$'' model).
While this problem arises in many practical applications and many techniques have been recently developed for its solution,
to the best of our knowledge nobody provided a rigorous theoretical analysis of this important setup. The purpose of the present paper is to
fill in this gap.  

We consider one of the most common settings, classification of high-dimensional normal vectors where, 
unlike standard assumptions, the number of classes could be  large.  
We derive  non-asymptotic  conditions   on effects of significant 
features, and the low and the upper bounds for distances between classes required for  successful feature selection 
and  classification with a given accuracy. Furthermore, we study an asymptotic setup
where the number of classes is diverging with the dimension of feature space and while 
the number of samples per class is possibly limited. 
We point out on an interesting and, at first glance, somewhat counter-intuitive phenomenon that a large number of classes 
may be a ``blessing'' rather than a ``curse'' since, in certain settings, 
the precision of classification can improve as the number of classes grows. This is due to more accurate feature selection 
since even weaker significant features, which are not sufficiently strong to be manifested in a coarse classification, being shared across the classes,
have a stronger impact as the number of classes increases.
We supplement our theoretical investigation by a simulation study and a real data example where we again observe the above phenomenon.  
\end{abstract}


\bigskip


\section{Introduction} \label{sec:intr}

Classification has been studied in many contexts.
In the era of ``Big Data'' one is usually interested in classifying objects that are described 
by a large number of features and belong to many different groups. 
For example the large hand-labeled ImageNet dataset {\tt http://www.image-net.org/} 
contains 10,000,000 labeled images depicting more than 10,000  object categories 
where each image, on the average, is represented by $482 \times 415 \approx 200,000$ pixels 
(see Russakovsky {\it et al.}, 2015 for description and discussion of this data set).
The challenge of handling large dimensional data got the name of ``large $p$ small $n$'' 
type of problems which means that dimensionality of parameter space $p$  by far exceeds the sample size $n$. 
It is well known that solving problems  of this type require rigorous model selection. 
In fact, the results of Bickel and Levina (2004), Fan and Fan (2008), Shao {\em et al.} (2011) 
demonstrate that even for the standard case of two classes, 
classification of high-dimensional normal vectors without feature selection 
is as bad as just pure random guessing.
However, while analysis of high-dimensional data (``Big data'') became ubiquitous, to the best of our knowledge, 
there are no theoretical studies that examine  the effect of large number of classes on classification accuracy.
The objective of the present paper is to fill in this gap.

At first glance, the problem of successful  classification when the number of classes is large seems close to impossible. 
On the other hand, humans have no difficulty in distinguishing between thousands of 
objects, and   the accuracy of state-of-the-art computer vision  techniques is approaching human accuracy.
In fact, in some settings, the accuracy of classification improves when the number of classes grows.
How is this possible? One of the reasons why multi-class classification succeeds is that 
selection of appropriate features from a large sparse $p$-dimensional vector becomes easier when 
the number of classes is growing since even weaker significant features that are not sufficiently strong to be manifested in a coarse
classification with a small number of classes may nevertheless have a strong impact
as the number of classes grows. Simulation studies in  Davis, Pensky and Crampton  (2011) 
and Parrish and Gupta (2012) support such a claim. Arias-Castro,  Cand\`{e}s  and  Plan (2011) 
reported on a similar occurrence for testing in the sparse  ANOVA model.   
Our paper  establishes a firm theoretical foundation under the above phenomenon 
and confirms it via simulation studies and a real data example.

Although there exists an enormous amount of literature on classification,
most of the existing theoretical results have been obtained for the binary classification ($L=2$) 
(see Boucheron, Bousquet and Lugosi, 2005 and references therein for a comprehensive survey). 
In particular, binary classification of high-dimensional sparse Gaussian vectors was considered 
in Bickel and Levina (2004), Fan and Fan (2008), Donoho  and  Jin  (2009 ab), 
Ingster,  Pouet  and  Tsybakov  (2009) and Shao {\em et al.} (2011) among others. 
 
In the  meantime,   a significant amount of effort has been spent on designing methods 
for the multi-class classification in statistical and machine learning literature. 
We can mention here techniques designed to adjust pairwise classification
to multi-class setting (Escalera {\it et al.}, 2011; Hill  and Doucet, 2007; Jain  and Kapoor, 2009), 
adjustment of the support vector machine technique to the case of several classes
(Crammer and Singer, 2001; Lee, Lin and Wahba, 2004) as well as a variety of approaches to expand the linear 
regression and the neural networks techniques to accommodate the multi-category setup  (see, e.g., Gupta, Bengio  and Weston, 2014). 
Tewari and Bartlett (2007) and Pan, Wang and Li (2016)   generalized theoretical  results for binary classification to the case of
multi-class classification and established consistency of the proposed classification procedures.
However, all above-mentioned investigations  considered only the ``small $L$, large $p$, small $n$'' setup, 
where the number of classes was assumed to be {\em fixed}.

This paper is probably the first attempt to rigorously investigate  ``large $L$, large $p$, small $n$'' classification and the impact of the number 
of classes on the accuracy of feature selection and classification.
In particular, we  explore the  somewhat  counter-intuitive phenomenon, where the large number of
classes may become a ``blessing'' rather than a ``curse'' for successful classification as more significant features may be revealed.
For this purpose,  we consider a well-known problem of multi-class classification of high-dimensional normal vectors.
We assume that only a subset of truly significant features really contribute to separation  between classes (sparsity). 
For this reason, we carry out feature selection and, following a standard scheme,
assign the new observed vector to the closest class w.r.t. the scaled Mahalanobis distance 
in the space of the selected significant features.
Our paper considers a realistic scenario where the number of classes as well as the number of features is large while 
the number of observations per class is possibly limited (``large $L$, large $p$, small $n$'' model).
We do not fix the total number of observations since in the real world the experience of each new class   
 comes with its own, usually finite, set of observations.

We start with a non-asymptotic setting and derive the conditions on effects of significant 
features, and the low and the upper bounds for the distances between classes required for  successful feature selection 
and  classification with a given accuracy.  All the results are obtained with the explicit constants and remain valid for any combination of parameters. 
Our finite sample study is followed by an asymptotic analysis for a large number of features $p$, 
where,  unlike previous works,   the number of classes $L$ may grow with $p$ while the number of samples per class 
may grow or stay fixed. Our findings indicate that having larger number of classes aids the feature selection and, hence, can improve 
classification accuracy. On the other hand,  larger  number of classes require  having larger  number of significant 
features $p_1$ for their separation which automatically leads to a ``large $p$'' setting. 
Nevertheless, due to increasing point isolation in high-dimensional spaces (see e.g. Giraud, 2015, Section 1.2.1),
those separation conditions become attainable when $p$ is large.

We ought to point out that our paper does not propose a novel
methodology for feature selection or classification. Rather than that, it studies one of the most popular Gaussian setting 
and adapts to the case of a large number of classes a standard general scheme, 
where feature selection is implemented by a thresholding technique
with the properly chosen threshold  and 
classification is carried out on the basis of the minimal Mahalanobis distance 
(we consider both the known and the unknown covariance matrix scenarios). 
This is a common widely used general scheme for classification and feature  selection in such setting
(see, e.g., Fan and Fan 2008; Shao {\em et al.}, 2011 and
Pan, Wang and Li, 2016 for similar approaches that differ mostly by selections of thresholds and distances).
Nevertheless, the setup is simple enough for derivations of  conditions required  for successful classification 
with a specified precision when the number of classes is large. 
Therefore, in our simulation study we do not compare these simple and well known techniques with the state of the art 
classification methodologies but instead investigate how these popular procedures perform when $p$ is large and
both the number of classes  $L$ and the number of significant features $p_1$  are growing. 
In particular, simulations support our finding that classification precision can improve when $L$ is increasing. 
The real data example confirms that the phenomenon above is not due to an artificial construction and is possible 
in a real life setting.

The rest of the paper is organized as follows.
In Section \ref{sec:model} we present the feature selection and multi-class
classification procedures and derive the non-asymptotic bounds for their accuracy. 
An asymptotic analysis is considered
in Section \ref{sec:asymp}.  Section \ref{sec:unknownS} 
discusses  adaptation of the procedure in the case of the unknown covariance matrix.  
In Section \ref{sec:examples}  we illustrate the
performance of the proposed approach on simulated and real-data examples.
Some concluding remarks are summarized in Section \ref{sec:remarks}.
All the proofs are given in the Appendix.


\section{Feature selection and classification procedure} \label{sec:model}

\subsection{Notation and preliminaries}
Consider the problem of multi-class classification of $p$-dimensional normal vectors with $L$ classes: 
\begin{equation} \label{eq:model0}
{\bf Y}_{li}={\bf m}_{l}+{\pmb \epsilon}_{li},\;\;\;\;l=1,\ldots, L; \;\;\;i=1,\ldots n_l,
\end{equation} 
where ${\bf m}_l \in \mathbb{R}^p$ is the vector of mean effects of $p$ features in the $l$-th class and
${\pmb \epsilon}_{li} \sim N({\bf 0}_p,\Sigma)$ with the common non-singular covariance matrix 
$\Sigma \in \mathbb{R}^{p \times p}$. To clarify the proposed approach we assume meanwhile that $\Sigma$ is known and discuss 
the situation with the unknown $\Sigma$ in Section \ref{sec:unknownS}.   

In what follows, we  study a realistic scenario where the number of classes as well as the number of features is large while 
the number of observations per class is possibly limited (``large $L$, large $p$, small $n$'' model).
We do not fix the total number of observations since in the real world the experience of each new class   
comes with its own, usually finite, set of observations.

After averaging over repeated observations within each class, model (\ref{eq:model0}) yields
\be \label{eq:model1}
\bar{\bf Y}_l={\bf m}_l+{\pmb \epsilon}^*_l,\;\;\;\;l=1,\ldots, L
\ee
where ${\pmb \epsilon}^*_l \sim N({\bf 0}_p, n_l^{-1}\Sigma)$.

The objective  is to assign a new observed  feature vector
${\bf Y}_0 \in \mathbb{R}^p$ to one of the $L$ classes. 
Denote
\begin{equation} \label{eq:rho}
N = \sum_{l=1}^L n_l, \quad  \rho_l = n_l/(n_l +1) \quad {\rm and} \quad L_1 = L-1,
\end{equation} 
where evidently  $1/2 \leq \rho_l < 1$.

Since $\Var({\bf Y}_0-\bar{\bf Y}_l)= \rho_l^{-1}\, \Sigma$,
we assign $\bf{Y}_0$ to the class $l$ with the nearest centroid 
$\bf{\bar{Y}}_l$  w.r.t to the scaled Mahalanobis distance:
\begin{equation} \label{eq:class0}
\hat{l}=\arg \min_{1 \leq l \leq L} \left\{\rho_l \, ({\bf Y}_0-\bar{\bf Y}_l)^t \Sigma^{-1}({\bf Y}_0-\bar{\bf Y}_l) \right\}.
\end{equation}

It is well-known (see, e.g., Bickel and Levina, 2004, Fan and Fan, 2008 and Shao {\em et al.}, 2011) that the performance
of classification procedures is worsening as the number of features grows (curse
of dimensionality). Hence, dimensionality reduction by feature selection prior  
to classification is crucial for large values of $p$.

Re-write (\ref{eq:model1}) in terms of the one-way multivariate analysis of variance (MANOVA) model as follows: 
\be \label{eq:model}
\bar{\bf Y}_l={\pmb \delta}+{\pmb \beta_l}+{\pmb \epsilon}^*_l,\;\;\;\;l=1,\ldots, L;
\ee
where ${\bf m}_l={\pmb \delta}+{\pmb \bbeta}_l$, ${\pmb \delta}$ is the vector of mean main effects of features
and $\beta_{lj},\;j=1,\ldots,p$ is the mean interaction effect of $j$-th feature with $l$-th class, with the standard identifiability conditions
$\sum_{l=1}^L \beta_{lj}=0$ for each $j=1,\ldots,p$.

The impact of $j$-th feature on classification depends 
on its variability between the different classes characterized by the interactions 
$\beta_{lj},\;l=1,\ldots,L$ in the model
(\ref{eq:model}). The larger are the interactions, the stronger is the impact of the feature.
A natural global measure of feature's contribution to classification is then
$b_j^2=\sum_{l=1}^L \beta^2_{lj}$. Note that a feature may still have a strong main effect 
$\delta_j$ but its contribution to classification nevertheless remains weak if it does not vary significantly
between classes, that is, if $b_j^2$ is small. 
The main goal of feature selection is to identify a sparse subset of significant features for 
further use in classification.


\subsection{Oracle classification} \label{subsec:oracle}

First, we consider an ideal situation where there is an oracle that provides the list of  truly significant 
features with $b_j^2>0$. In this case, we would obviously use only those 
features for classification, thus, reducing the dimensionality of the problem.
Define indicator variables $x_j=I\{b_j^2 > 0\}$, and
let $p_1=\sum_{j=1}^p x_j$ and $p_0=p-p_1$ be,  respectively, the numbers of significant and non-significant 
features. Without loss of generality, we can always order features in such a way
that those $p_1$ significant features are the first ones. The classification procedure (\ref{eq:class0})
then  becomes  
\begin{equation} \label{eq:trueclass}
\hat{l}= \underset{1 \leq l \leq L}{\operatorname{argmin}}
\left\{\rho_l \, ({\bf Y}^*_0-\bar{\bf Y}^*_l)^t (\Sigma^*)^{-1}({\bf Y}^*_0-\bar{\bf Y}^*_l) \right\},
\end{equation}
where ${\bf Y}_0^*, {\bf Y}_l^* \in \mathbb{R}^{p_1}$ are the truncated versions
of ${\bf Y}_0$ and $\bar{\bf Y}_l$ respectively:  $Y_{0j}^*=Y_{0j}$ and
$Y_{lj}^*=\bar{Y}_{lj},\;j=1,\ldots,p_1$, and $\Sigma^* \in \mathbb{R}^{p_1 \times p_1}$ 
is the corresponding upper left sub-matrix of $\Sigma$.
 
Theorem \ref{th:oracle} provides an upper bound for misclassification 
error of the oracle classification procedure (\ref{eq:trueclass}):

\begin{theorem} \label{th:oracle}
Consider the model (\ref{eq:model0}) and the equivalent model (\ref{eq:model}). 
Let ${\bf m}^*_k \in \mathbb{R}^{p_1},\;k=1,\ldots,L,$ be the truncated versions of 
class centers ${\bf m}_k$ and assume that for all pairs of classes 
\be \label{eq:as1}
({\bf m}^*_k-{\bf m}^*_{k'})^t (\Sigma^*)^{-1}({\bf m}^*_k-{\bf m}^*_{k'}) \geq \frac{8  \, \ln(L_1/ \alpha)}
{\min(\rho_k,\rho_{k'})}\
\lkv 1 + \frac{1}{\sqrt{2 \min(n_k, n_{k'})}} \lkr 1 + \sqrt{\frac{2 p_1}{\ln(L_1/ \alpha)}}\rkr \rkv
\ee
for some $0 < \alpha \leq 1 $. 

Let a new observation $Y_0$ from  the class $l$ be assigned  to 
the $\hat{l}$-th class according to classification rule (\ref{eq:trueclass}). Then,
the misclassification error is
\be \label{eq:misclass1}
P(\hat{l} \neq l) \leq \alpha
\ee
\end{theorem}

Condition (\ref{eq:as1}) verifies that classes should be sufficiently separated from each other 
(in terms of Mahalanobis distance) to achieve the required classification accuracy. In fact, the requirements in (\ref{eq:as1})  are also essentially necessary.
Theorem \ref{th:lower1} below, which is a direct consequence of Fano's lemma  
for the lower bound of misclassification error (see, e.g., Ibragimov and Hasminskii, 1981, Section 7.1),  
implies that the first term $O\left(\ln(L_1/\alpha)\right)$ in the RHS of \fr{eq:as1} is  unavoidable for successful
classification and cannot be significantly improved (in the minimax sense) even in the idealized case, where the class centers ${\bf m}^*_k$ are known:
\begin{theorem} \label{th:lower1}
Consider the model (\ref{eq:model0}). Let a new
observation $\pmb{Y}_0$ be from one of $L$ classes. If  
\be \label{eq:lower}
\tilde{\Delta}^2 =  \min_{l \ne k}~ ({\pmb m}^*_l - {\pmb m}^*_k)^t (\Sigma^*)^{-1}
({\pmb m}^*_l - {\pmb m}^*_k)
\leq    2 \, \aleph    \ln L_1
\ee
for some $\aleph > 0$,  then
\be \label{eq:lowbou}
\inf_\psi \max_{1 \leq l \leq L} P_l(\psi(\pmb{Y}_0) \neq l) \geq 
1-\aleph -\frac{\ln 2}{\ln L_1},
\ee
where $P_l$ is the probability evaluated under the assumption that $\pmb{Y}_0$ belongs to the $l$-th class, 
and  the infimum is taken over all classification
rules $\psi(\pmb{Y}_0): \pmb{Y}_0 \rightarrow \{1,\ldots,L\}$. 
\end{theorem}

The second term  in the RHS of (\ref{eq:as1}) appears due to replacing the unknown $p_1$-dimensional
class centers ${\bf m}^*_k$'s by the corresponding within-class sample means $\bar{\bf Y}^*_k$'s in (\ref{eq:trueclass}). 
Indeed,  straightforward extension of  the results of Theorem 1 of Fan and Fan (2008) for a general $L \geq 2$ yields that, unless for all 
pairs $(k,k')$, $({\bf m}^*_k-{\bf m}^*_{k'})^t (\Sigma^*)^{-1}({\bf m}^*_k-{\bf m}^*_{k'}) \geq C \sqrt{\frac{p_1 \ln L_1}{\min(n_k,n_{k'})}}$ for some
$C>0$, 
the curse of dimensionality affects the accumulated error in estimating high-dimensional ${\bf m}^*_k$'s
and yields classification performance nearly the same as random guessing.


\subsection{Feature selection procedure} \label{subsec:selection}

Consider now classification setup in the MANOVA model (\ref{eq:model}) with a more realistic scenario, 
where a set of significant features is unknown and should be identified from the data.

To simplify the calculus and to avoid complications with post-selection inference, we split the data at random
into two sets $Y_{lj}^{(1)}$'s and $Y_{lj}^{(2)}$'s in some fixed proportion $\phi \in (0,1)$ 
(in the simplest case, the sizes of both sets are equal with $\phi=1/2$). Subsequently,  use $Y_{lj}^{(1)}$'s
for feature selection and $Y_{lj}^{(2)}$'s for classification based on the selected features.
More specifically, for $l$-th class, split its $n_l$ observations $Y_{lj}$'s into two sub-samples of sizes 
$n^{(1)}_l$ and $n^{(2)}_l$ at the same proportion $\pi$, i.e.
$n^{(1)}_l=\lfloor \pi n_l \rfloor$, where $\lfloor \cdot \rfloor$ is the integer part,  
and $n^{(2)}_l=n_l-n^{(1)}_l,\;l=1,\ldots,L$.
Denote the total sample sizes of the resulting two sets
by $N_1=\sum_{l=1}^L n^{(1)}_l$ and $N_2=\sum_{l=1}^L n^{(2)}_l$, 
so that $N_1 + N_2 = N$.  

Following our previous arguments, 
a $j$-th feature is not significant (irrelevant) for classification if it has zero
interaction effects with all classes, that is, if
$\beta_{lj}=0,\;j=1,\ldots,L$ or, equivalently, $b_j^2=0$.
Then, for  each $j=1,\ldots,p$ we need  to test the null hypothesis
$H_{0j}:b_j^2=0$.
An obvious test statistic is then 
\be \label{eq:zetaj}
\zeta_j = \sigma_j^{-2}~ \sum_{l=1}^L n_l^{(1)} (\bar{Y}^{(1)}_{lj}-\bar{Y}^{(1)}_{\cdot j})^2,
\ee
where $\sigma^2_j=\Sigma_{jj}$ and $\bar{Y}^{(1)}_{\cdot j}= (n_l^{(1)})^{-1}\, \sum_{l=1}^L Y^{(1)}_{lj}$.
Under the null, $\zeta_j \sim \chi^2_{L_1}$, while under the
alternative $\zeta_j \sim  \chi^2_{L_1;\mu_j}$, where
$\chi^2_{L_1;\mu_j}$ is the non-central chi-square distribution with the
non-centrality parameter $\mu_j=\sigma_j^{-2}~\sum_{l=1}^L n^{(1)}_l \beta_{lj}^2$. Note that unless $\Sigma$ is diagonal, $\zeta_j$'s are correlated.

For a given $0 < \alpha \leq 1$, define a threshold
\be \label{eq:lambda}
\lambda = L_1+2\sqrt{L_1 \ln (2p/\alpha)} +2\ln (2p/\alpha)
\ee
and select the $j$-th feature as significant
(reject $H_{0j}$) if 
\begin{equation} \label{eq:test}
\zeta_j = \sigma^{-2}_j~\sum_{l=1}^L n^{(1)}_l (\bar{Y}^{(1)}_{lj}-\bar{Y}^{(1)}_{\cdot j})^2 > \lambda 
\end{equation}

The following theorem shows that under certain conditions on the minimal
required effect for significant features, the proposed feature selection procedure
correctly identifies the true (unknown) subset of significant features with probability at least $1-\alpha$:

\begin{theorem} \label{th:pi1}
Consider the feature selection procedure (\ref{eq:test}) with the threshold (\ref{eq:lambda}) for some $0 < \alpha \leq 1$. 
Define indicator variables $\hat{x}_j=I\{\sigma^{-2}_j~\sum_{l=1}^L n^{(1)}_l (\bar{Y}_{lj}-\bar{Y}^{(1)}_{\cdot j})^2 > \lambda\},\;j=1,\ldots,p$. 
Let  
\be \label{eq:mu*}
\mu^* = \min_{1 \leq j\leq p_1} \sigma^{-2}_j~\sum_{l=1}^L n^{(1)}_l \beta_{lj}^2
\ee
and assume that for all $p_1$ truly significant features one has
\be \label{eq:as2}
\mu^*    \geq 4 \left(3 \ln (2p/\alpha) + \sqrt{L_1 \ln (2p/\alpha)}\right) 
\ee
Then, 
$$
P(\hat{x}=x) \geq 1 - \alpha
$$
\end{theorem}

The condition (\ref{eq:as2}) on the total minimal effect for significant features
can be re-formulated in terms on their {\em average} effect per class:
\be \label{eq:per_class}
\frac{1}{\sigma^2_j L}~\sum_{l=1}^L n^{(1)}_l \beta_{lj}^2  \geq 4 \left(\frac{3 \ln (2p/\alpha)}{L} + 
\sqrt{\frac{\ln (2p/\alpha)}{L}}\right),\;\;\;j=1,\ldots,p_1
\ee
Thus, as the number of classes in model (\ref{eq:model0}) increases, even significant features with weaker
effects within each class become manifested and contribute to classification. 
Effect of a certain feature that remains latent and unnoticed in coarse 
classification with a small number of classes may be expressed in a finer
classification.


\subsection{Classification rule and misclassification error} \label{subsec:class}
Consider now the classification rule (\ref{eq:trueclass}) applied on the second set of the data with $\bar{Y}_l^{(2)*}$, where the unknown true
$x_j$ are replaced by $\hat{x}_j$ following the proposed feature selection procedure.
Let
$\hat{p}_1=\sum_{j=1}^p \hat{x}_j$ be the number of features declared significant and $\hat{p}_0=p-\hat{p}_1$. 
Again, order the features in such a way that those $\hat{p}_1$ features selected as significant 
are the first ones. Thus, the  resulting classification rule can then be presented as follows:
\begin{equation} \label{eq:class}
\hat{l}= \underset{1 \leq l \leq L}{\operatorname{argmin}} 
\left\{\rho_l \, ({\bf Y}^*_0-\bar{\bf Y}^{(2)*}_l)^t (\Sigma^*)^{-1}({\bf Y}^*_0-\bar{\bf Y}^{(2)*}_l) \right\},
\end{equation}
where the truncated vectors ${\bf Y}^*_0, \bar{\bf Y}^{(2)*}_l \in \mathbb{R}^{\hat{p}_1},\;l=1,\ldots,L$ 
are defined now as $Y^*_{0j}=Y_{0j},\; Y^{(2)*}_{lj}=\bar{Y}^{(2)}_{lj},\;j=1,\ldots,\hat{p}_1$, and
$\Sigma^* \in \mathbb{R}^{{\hat p}_1 \times {\hat p}_1}$ is the corresponding upper left sub-matrix of
$\Sigma$, and $\rho_l=n_l^{(2)}/(n_l^{(2)}+1)$. 

We have
\be \label{eq:error}
P(\hat{l} \neq l) \leq P(\hat{l} \neq l \mid \hat{x}=x) + P(\hat{x} \neq x),
\ee
where, due to the fact that   different data was used for feature selection and classification, 
by Theorem \ref{th:oracle} and Theorem \ref{th:pi1}, each probability
in the RHS of (\ref{eq:error}) is at most $\alpha$. Thus, the following result holds:
\begin{theorem} \label{th:main}
Consider the model (\ref{eq:model0}) and the corresponding model (\ref{eq:model}).
Assume the conditions (\ref{eq:as1}) (with $n_l$ replaced by $n_l^{(2)}$) and (\ref{eq:as2}) hold for some
$0 < \alpha \leq 1/2$. 
Apply feature selection procedure
(\ref{eq:test}) and use the selected features for classification via the rule 
(\ref{eq:class}). Then, 
$$
P({\rm correct\ classification}) \geq 1 - 2\alpha
$$
\end{theorem}


\section{Asymptotic analysis}  \label{sec:asymp}

Conditions \fr{eq:as1} and \fr{eq:as2} (or \fr{eq:per_class}) of 
Theorems 1 and 2, respectively, provide the  non-asymptotic  lower bounds 
on the minimal distance  between different  classes and the minimal effect of   significant 
features required for the perfect feature selection and classification error bounded above by $2 \alpha$. 
In order to gain better understanding of these conditions, we consider an
asymptotic setup. 

Standard asymptotics considered in classification literature assume 
that the number of features $p$ and the sample sizes $n_l$ increase whereas the
number of classes $L$ is fixed  (see,   e.g., Fan and Fan, 2008; 
Shao {\em et al.}, 2011 for $L=2$ and Pan, Wang and Li, 2016, for a general but fixed $L$).
On the contrary, our study is motivated by the case where the number of classes may
also be large (``large $L$, large $p$, small $n$'').

Recall that $N=\sum_{l=1}^L n_l$  is the total sample size and let
the number of features $p \to \infty$.  Following Pan, Wang and Li (2016), assume that all eigenvalues 
of the $p_1 \times p_1$ covariance matrix of significant features $\Sigma^*$ are finite and bounded away from zero, i.e.,
there exist absolute constants $\tau_1$ and $\tau_2$ such that
\be \label{eq:eigenassump}
0 < \tau_1 \leq  \lambda_{\min} (\Sigma^*) \leq  \lambda_{\max} (\Sigma^*) \leq \tau_2 < \infty.
\ee
The samples sizes $n_l$ within classes also grow with $p$. For simplicity of 
exposition, we assume that they are of the same asymptotic order and splitted more or less equally between the two sets ($\pi \sim 1/2$), that is, $n^{(1)}_l \sim n^{(2)}_l \sim n$ for all $l=1,\ldots,L$, where  
$n=N/(2L)$ and $a \sim b$ means $a=b(1+o(1))$. In such asymptotic setup, 
$\rho_l \sim 1 - 1/n$, while $\sqrt{1-\rho_l \rho_k} \sim \sqrt{2/n}$. 
Though the results in the previous section allow one to study various other settings with unequal class sizes,  the asymptotic analysis of
a vast variety of such possible scenarios is beyond the scope of this paper.
 

Consider now the condition (\ref{eq:as1}) of Theorems \ref{th:oracle} and \ref{th:main} on the minimal separation Mahalanobis distance between
any two class centers as $p$ tends to infinity, while $n$, the number
of significant features $p_1$ and the number of classes $L$ may increase with $p$,
and $\alpha$ may depend on $n, p$ and $L$. Thus, (\ref{eq:as1}) yields:
\be \label{eq:as1asymp}
\min_{k \neq k'}~ ({\pmb m}^*_k-{\pmb m}^*_{k'})^t (\Sigma^*)^{-1} ({\pmb m}^*_k-{\pmb m}^*_{k'}) \geq \Delta^2_* \sim 8 \ln(L_1/\alpha)
\left(1+\frac{1}{\sqrt{2n}}\left(1+\sqrt{\frac{2p_1}{\ln(L_1/\alpha)}}\right)\right)
\ee
%
Define
$$
\eta_1 = \lim_{p \to \infty} \sqrt{\frac{p_1}{n \ln(L_1/\alpha)}}
$$
Depending on $\eta_1$, the condition (\ref{eq:as1asymp}) implies two possible asymptotic regimes for $\Delta^2_*$: 
 \be \label{eq:asymp_distance}
\Delta^2_* \sim  \lfi \begin{array}{lll}
8 \ln \lkr \frac{L_1}{\alpha}\rkr(1 + \eta_1),\! & \! 0 \leq  \eta_1 < \infty \! & \! \mbox{(sparse regime - small number  of  significant  features)}\\
8 \, \sqrt{\frac{p_1 \ln(L_1/\alpha)}{n}},  \! & \!\eta_1= \infty 
\! & \! \mbox{(dense  regime - large number  of  significant  features)} 
\end{array} \right.
\ee
For sparse regime
($\eta_1 < \infty$), the required minimal between-class distance $\Delta_*^2$ grows slowly as $\ln L$ and
from Theorem \ref{th:lower1} it immediately follows that 
this is the lowest possible rate
for successful classification:  

\begin{proposition} \label{prop:lower1}
Let $L  \to \infty$ and $p_1 \to \infty$ as $p \to \infty$. 
Let a new observation $\pmb{Y}_0$ be from one of $L$ classes. If
$$
\Delta^2_* \sim  2\, \delta_{p_1}    \ln L_1,
$$
where $\delta_{p_1} \to 0$ arbitrarily slow as $p \to \infty$, then 
$$
\lim_{p \to \infty} \inf_\psi \max_{1 \leq l \leq L} P_l(\psi(\pmb{Y}_0) \neq l) = 1, 
$$
where $P_l$ is the probability evaluated under the assumption that $\pmb{Y}_0$ belongs to the $l$-th class, 
and  the infimum is taken over all classification
rules $\psi(\pmb{Y}_0): \pmb{Y}_0 \rightarrow \{1,\ldots,L\}$. 
\end{proposition}

For dense regime, the number of significant features $p_1$ is large
enough for the accumulated error of estimating $p_1$-dimensional 
${\bf m}^*_k$'s by $\bar{\bf Y}^{(1)*}_k$'s to become
dominant (see Section \ref{subsec:oracle}) and the classes should be, 
therefore, much stronger separated to deal with the curse of dimensionality.

It is natural that for successful classification the between-class distances should grow with $L$. Note, however, that unless the
number of classes increases exponentially with $p_1$, the growth rate of $\Delta_*^2$ is $o(p_1)$
and the corresponding average per-feature distances 
$\frac{1}{p_1}({\pmb m}^*_k-{\pmb m}^*_{k'})^t (\Sigma^*)^{-1}({\pmb m}^*_k-{\pmb m}^*_{k'}) $ still tend to zero.

Similarly, from the condition (\ref{eq:as2}) in Theorems~\ref{th:pi1}~and~\ref{th:main} 
on the minimal effect for significant features required  for the
perfect feature selection, we have asymptotically
$$
b^2_*=\min_{1 \leq j \leq p_1} \sigma^{-2}_j b_j^2 \sim \frac{4}{n}
\left(3 \ln(2p/\alpha)+\sqrt{L_1 \ln(2p/\alpha)}\right)
$$
%
Let
$$
\eta_2 = \lim_{p \to \infty} \sqrt{\frac{\ln(2p/\alpha)}{L_1}}
$$
Then, 
\be \label{eq:as2asymp}
b_*^2 \sim \lfi \begin{array}{lll}
4 n^{-1}\, \sqrt{L_1 \ln(2p/\alpha)}(1 + 3 \eta_2), & 0 \leq  \eta_2 < \infty
& ({\rm large\;number\;of\;classes})\\ 
12  n^{-1}\, \ln(2p/\alpha), & \eta_2 = \infty 
& ({\rm small\;number\;of\;classes})
\end{array} \right.
\ee 
and the threshold $\lambda$ in (\ref{eq:lambda}) for feature selection can be presented as
$$
\lambda \sim \lfi \begin{array}{ll}
L_1 (1+2\eta_2+2\eta_2^2), & 0 \leq \eta_2 < \infty \\
2\ln(2p/\alpha), & \eta_2=\infty
\end{array} \right.
$$

To gain some insight on the minimal required effect for a significant feature 
to contribute to classification as the number of classes increases, assume for simplicity that each significant feature has 
equal effects on each class, that is, $\beta_{lj}$ in (\ref{eq:model})   vary only in signs:  
$\beta^2_{lj}= \beta^2_j$, $l=1, \ldots, L$.  
%
Since $0 \leq  \eta_2 < \infty$ implies that $L$ is large, so that $L_1=L-1 \sim L$, condition (\ref{eq:as2asymp}) yields as $p \to \infty$: 
\be \label{eq:as2asymp1}
\beta_j^2 \sim \lfi \begin{array}{lll}
4 \sigma_j^2~ n^{-1}\, \eta_2 (1 + 3 \eta_2), & 0 \leq  \eta_2 < \infty
& ({\rm large\;number\;of\;classes})\\ 
12 \sigma_j^2~ n^{-1} \, L^{-1}  \ln(2p/\alpha), & \eta_2 = \infty 
& ({\rm small\;number\;of\;classes})
\end{array} \right.
\ee
Since $\eta_2$ is decreasing with $L$ for a given value of $\alpha$, the required minimal level for $\beta_j^2$ in the RHS of (\ref{eq:as2asymp1}) 
decreases as $L$ grows and, therefore, more significant
features   become manifested in classification for larger number of classes. 
Thus, while it might be hard to perform coarse classification 
with a set of weak features, their impacts grow  as one considers finer and finer separation between objects 
(see also the corresponding remarks at the end of Section \ref{subsec:selection}).

Although in this section our goal was to explore the case when $L \to \infty$, 
calculations above remain valid for a fixed value of $L$ (commonly,  $L=2$).
In particular, if $L$ is fixed and  $n=o(p)$,   conditions (\ref{eq:as1asymp}) and (\ref{eq:as2asymp1}) 
are of the form $\Delta^2_* \sim C_1\sqrt{\frac{p_1}{n}}$ and
$\beta^2_j \sim C_2 n^{-1}\ln(p/\alpha),\;C_1, C_2 > 0$ and are similar to those of
Fan and Fan (2008, Theorem 1 and Theorem 3). See also the results of
Donoho and Jin (2009 a,b) and
Ingster, Pouet and Tsybakov (2009) for closely related setups.


\section{Unknown covariance matrix}  
\label{sec:unknownS}

So far the covariance matrix $\Sigma$ was assumed to be known. In practice, however, it
should usually be estimated from the data. 
The standard MLE estimator 
based on the first sub-sample
\be \label{eq:sigmamle}
\widehat{\Sigma}^{(1)}=\frac{1}{N_1}\sum_{l=1}^L \sum_{i=1}^{n_l^{(1)}} \lkr {\bf Y}_{il}^{(1)} -\bar{\bf Y}_l^{(1)} \rkr
\lkr {\bf Y}_{il}^{(1)} -\bar{\bf Y}_l^{(1)} \rkr^t
\ee
and the similar unbiased pooled estimator commonly used in MANOVA behave poorly for high-dimensional data. 
However, under the sparsity assumption, the proposed classification procedure requires only to estimate
the variances $\sigma^2_j$ in  feature selection procedure  (\ref{eq:zetaj}) and the inverse of the upper left 
sub-matrix $\Sigma^* \in \mathbb{R}^{\hat{p}_1 \times \hat{p}_1}$ of $\Sigma$ in   classification rule  (\ref{eq:class}). 
Thus, when $p_1 \ll p$, a low-dimensional matrix $(\widehat{\Sigma^*})^{-1}$ may still be a
good estimator of the true sub-matrix $(\Sigma^*)^{-1}$ and (under some additional  mild conditions) 
may be used instead of the latter in (\ref{eq:class}).

Assume that $p \leq \frac{\alpha}{2}~e^{(N_1-L)/4}$. 
Replace  $\sigma^2_j$ in (\ref{eq:zetaj}) by $\hat{\sigma}^2_j=\widehat{\Sigma}^{(1)}_{jj}$ 
and consider the feature selection procedure (\ref{eq:test}) with a somewhat larger threshold 
\be \label{eq:lambda1}
\lambda_1=\frac{\lambda}{1-\kappa},
\ee
where $\lambda$ is the threshold (\ref{eq:lambda}) used for the case
of known variances and
\be \label{eq:kappa}
\kappa=\kappa(p,N_1,L,\alpha)=2\sqrt{\frac{\ln(2p/\alpha)}{N_1-L}}+2~\frac{\ln(2p/\alpha)}{N_1-L} < 1
\ee 
The following theorem  shows that under slightly stronger 
conditions on the minimal required effect for significant features, the above feature selection procedure with
estimated $\sigma_j^2$ still controls the probability of correct identification
of the true subset of significant features.

\begin{theorem} \label{th:pi2} 
Let $0 < \alpha \leq 1/2$   and assume that $p \leq \frac{\alpha}{2}~e^{(N_1-L)/4}$. 
Define indicator variables 
\be \label{eq:test1}
\hat{x}_j=I\{\hat{\sigma}^{-2}_j~\sum_{l=1}^L n^{(1)}_l (\bar{Y}^{(1)}_{lj}-\bar{Y}^{(1)}_{\cdot j})^2 > \lambda_1\},\;j=1,\ldots,p
\ee
with $\lambda_1$ given in (\ref{eq:lambda1}). 
%
Assume that $\mu_*$ in \eqref{eq:mu*} satisfies
\be \label{eq:as22}
\mu_*+L_1-2\sqrt{(L_1+2\mu_*)\ln(2p/\alpha)} > \lambda_1 (1+\kappa)
\ee
Then, 
$$
P(\hat{x}=x) \geq 1 - 2\alpha
$$
\end{theorem}

Consider now the classification procedure (\ref{eq:class}). In what follows we assume that $\Sigma^*$ is non-singular.
Consider an estimator $\widehat{\Sigma^*}$ of $\Sigma^*$ of the form
\be \label{eq:sigmamle0}
\widehat{\Sigma^*}=\frac{1}{N_2}\sum_{l=1}^L \sum_{i=1}^{n^{(2)}_l} ({\bf Y}^{(2)*}_{il}-\bar{\bf Y}^{(2)*}_l)
({\bf Y}^{(2)*}_{il}-\bar{\bf Y}^{(2)*}_l)^t,
\ee 
where ${\bf Y}^{(2)*}_{il}$ are the corresponding $\hat{p}_1$-dimensional truncated
versions of ${\bf Y}^{(2)}_{il}$.

Assign ${\bf Y}_0$ the $\hat{l}$-th class by replacing the true (unknown)
$(\Sigma^*)^{-1}$ in (\ref{eq:class}) by $(\widehat{\Sigma^*})^{-1}$:
\begin{equation} \label{eq:class1}
\hat{l}= \underset{1 \leq l \leq L}{\operatorname{argmin}} 
\left\{\rho_l \, ({\bf Y}^*_0-\bar{\bf Y}^{(2)*}_l)^t (\widehat{\Sigma^*})^{-1}({\bf Y}^*_0-\bar{\bf Y}^{(2)*}_l) \right\}.
\end{equation}
Then   the following version of Theorem \ref{th:main}  holds:
\begin{theorem} \label{th:main1}
Consider the model (\ref{eq:model0}) and the corresponding model (\ref{eq:model}), where  
$p \leq \frac{\alpha}{2}~e^{(N_1-L)/4}$, 
\be \label{eq:p1N}
\max\left(L,2\ln\left(\frac{2}{\alpha}\right)\right) < p_1 < 
\frac{1}{4C_1}\left(\frac{\lambda_{\min}(\Sigma^*)}{\lambda_{\max}(\Sigma^*)}\right)^4 N_2
\ee 
for some $0 < \alpha < 1/4$ and $C_1$ is an  absolute constant  specified  in the proof. 
Denote 
\be \label{eq:betap1n}
\gamma_{p_1,N_2}=2\frac{\lambda^2_{\max}(\Sigma^*)}{\lambda^2_{\min}(\Sigma^*)}\sqrt{\frac{C_1 p_1}{N_2}}
\ee
and note that $\gamma_{p_1,N_2}  < 1$ due to (\ref{eq:p1N}).
Assume the condition (\ref{eq:as22}) and a somewhat stronger version of the condition (\ref{eq:as1}), namely,
\be \label{eq:as11}
\begin{split}
({\bf m}^*_k-{\bf m}^*_{k'})^t (\Sigma^*)^{-1}({\bf m}^*_k-{\bf m}^*_{k'}) & \geq \frac{8  \, \ln(L_1/ \alpha)}
{(1-\gamma_{p_1,N_2})\min(\rho_k,\rho_{k'})} \\
& \times \lkv 1 + \sqrt{\frac{1}{2 \min \lkr n^{(2)}_k, n^{(2)}_{k'}\rkr }+ \gamma^2_{p_1,N_2}} \cdot
 \lkr 1 + \sqrt{\frac{2 p_1}{\ln(L_1/ \alpha)}}\rkr \rkv
\end{split}
\ee
Apply feature selection procedure
(\ref{eq:test1}) and use the selected features for classification via the rule 
(\ref{eq:class1}). Then,
$$
P({\rm correct\ classification}) \geq 1 - 4\alpha
$$
\end{theorem}

Theorem \ref{th:main1} shows that for a sparse setup the
proposed classification procedure can still be used when the covariance matrix
is unknown and estimated from the data.


\section{Examples} \label{sec:examples}
In this section we demonstrate the performance of the proposed feature selection
and classification procedure on simulated and real-data examples.
Its main goal is to illustrate the phenomenon of improving the accuracy as the number of classes grows discussed in the previous 
sections. 

We found that in practice there is no real need
to split the original data  and used the entire data set for
both feature selection and classification.

\subsection{Simulation study}  \label{subsec:simul}
Simulated examples follow the settings  presented in Pan, Wang and Li (2016).
 
We generated the class means  as i.i.d. normal vectors $\pmb{m}_l \sim N(0, \sigma_m^2 X),\;l=1,\ldots,L$, 
where $ X_{p \times p}$ is a diagonal matrix with $x_i =1$ for   $p_1$ indices and $x_i = 0$ for others. 
Since the vectors generated in this manner do not necessarily satisfy our assumptions, in order to reduce an impact of 
a particular choice of vectors $\pmb{m}_l$,  we generated $M_1$ replications of the class means. 
Furthermore, following the model (\ref{eq:model1}), for each replication of class means
$\pmb{m}_l,\;l=1,\ldots,L$  we generated $M_2$ sets of training samples 
$\bar{Y}_{lji} = m_{lj} + \epsilon^*_{lji},\;j=1, \ldots, p;\;i=1,\ldots, n$, 
where $\epsilon^*_{lji}$ are i.i.d. $N(0, n^{-1} \Sigma)$. Finally, for each of 
$M_1 \cdot M_2$ sets of training samples, we drew a test set of $M_3$ new vectors 
from randomly chosen classes as i.i.d. normal vectors $N(\pmb{m}_l, \Sigma)$.

   \begin{figure}
    \begin{center} 
    \includegraphics[scale=.7]{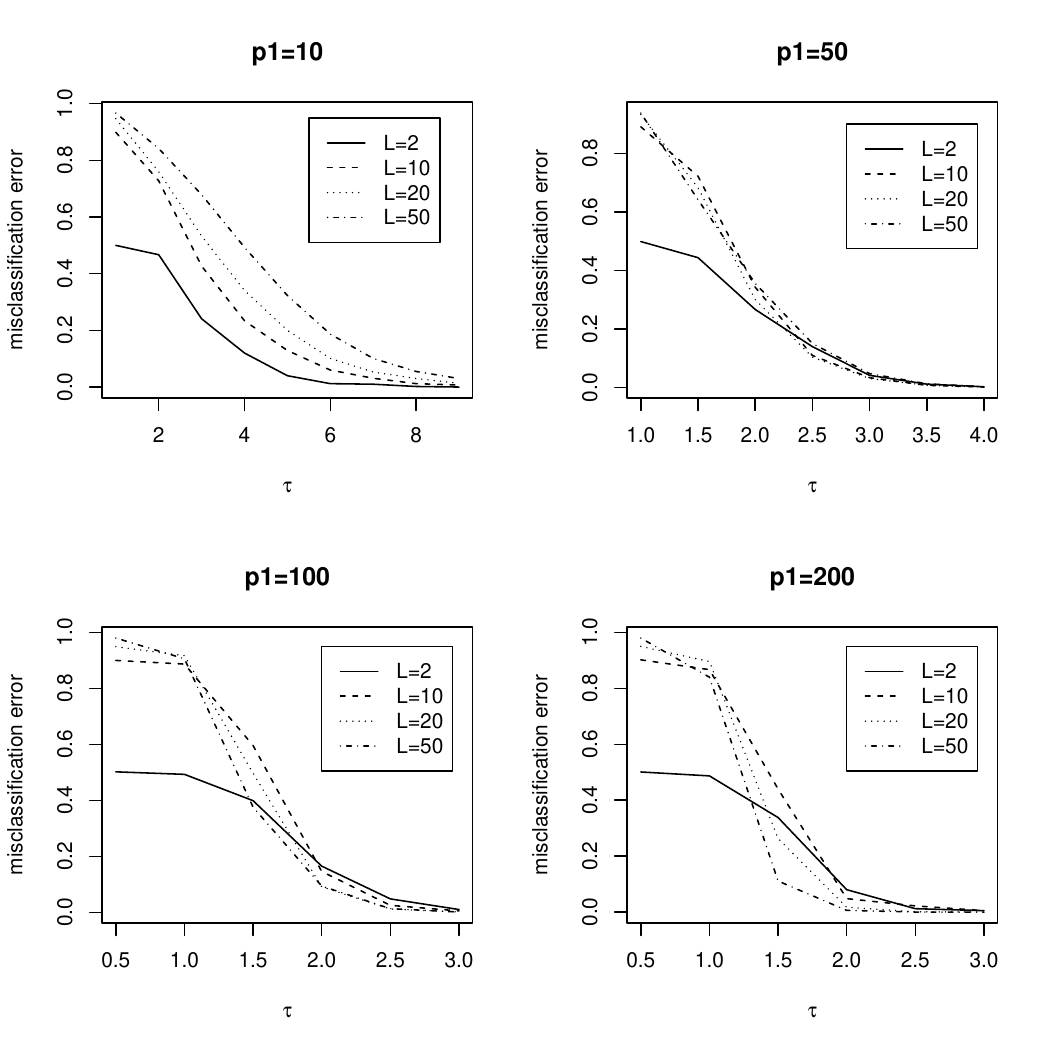}
   \caption{Average misclassification errors as  functions of $\tau$ for various combinations of $p_1$ and $L$ for Example 1.
  \label{fig:errors1} }
   \end{center}
   \end{figure}

 \begin{figure}[h] 
    \begin{center}
    \includegraphics[scale=.7]{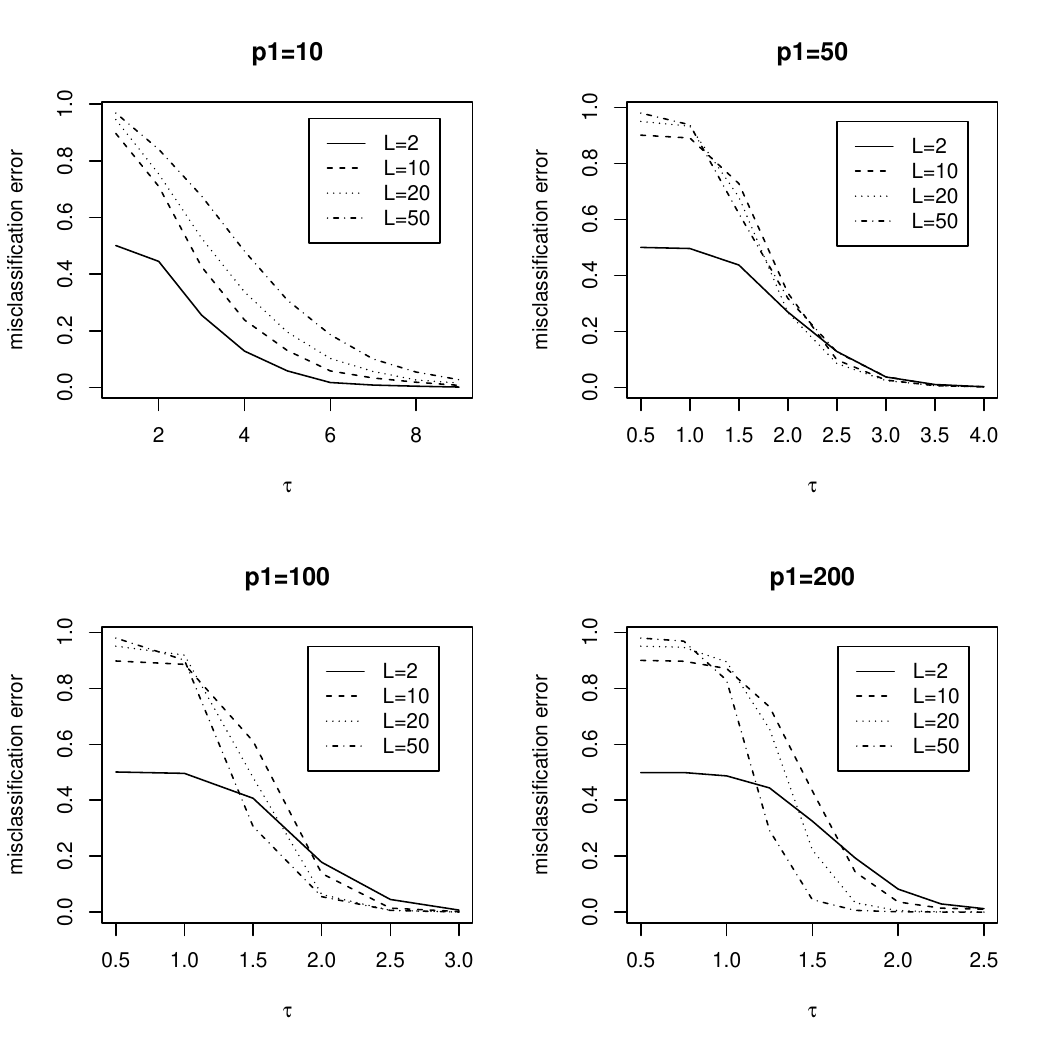}
   \caption{Average misclassification errors as  functions of $\tau$ for various combinations of $p_1$ and $L$ for Example 2.
  \label{fig:errors2} }
   \end{center}
   \end{figure}

We used the same three choices for covariance matrix $\Sigma$ as in Pan, Wang and Li (2016).
In Example 1 features were independent, i.e. $\Sigma=\sigma^2 I_p$. In Example 2 we used the 
autoregressive covariance structure with $\Sigma_{h_1,h_2}=\sigma^2~0.5^{|h_1-h_2|}$, while 
in Example 3 we set  $\Sigma_{h_1,h_2}=\sigma^2~(0.5+0.5 I\{h_1=h_2\}),\;h_1, h_2=1,\ldots,p $ implying equal
variances $\sigma^2$ and all covariances equal to $\sigma^2/2$ (compound symmetric structure). 
We carried out  simulations with both the true covariance matrix $\Sigma$ and its MLE $\widehat{\Sigma}$ given by (\ref{eq:sigmamle}). 
Since the performances of feature selection and classification procedures in both cases were similar, in  
what follows we present only the results obtained with   $\widehat{\Sigma}$.
 
For each training sample we first carried out the feature selection procedure described above 
with the threshold $\lambda_1$ defined in \fr{eq:lambda1} and $\alpha = 0.05$. Subsequently, we used the selected features 
for classifying $M_3$ vectors from the corresponding test set according to the rule \fr{eq:class1}. 
In the case when it delivered a non-unique solution, we chose one of the suggested solutions at random.

In all simulations we used $M_1 = M_2 = M_3 = 50,\; p=500,\; \sigma=1$ and $n=20$. 
Note that classification precision depends on the variance ratio $\tau^2 = \sigma_m^2/(\sigma^2/n)$ 
that may be viewed as a signal-to-noise ratio. 
For this reason, we studied performance of feature selection and classification 
for various combinations of $p_1$, $L$ and $\tau$. In particular, we used
$p_1=10, 50, 100, 200$, $L=2, 10, 20, 50$ and several values of $\tau$ depending on $p_1$.

\begin{figure}
    \begin{center}
    \includegraphics[scale=.7]{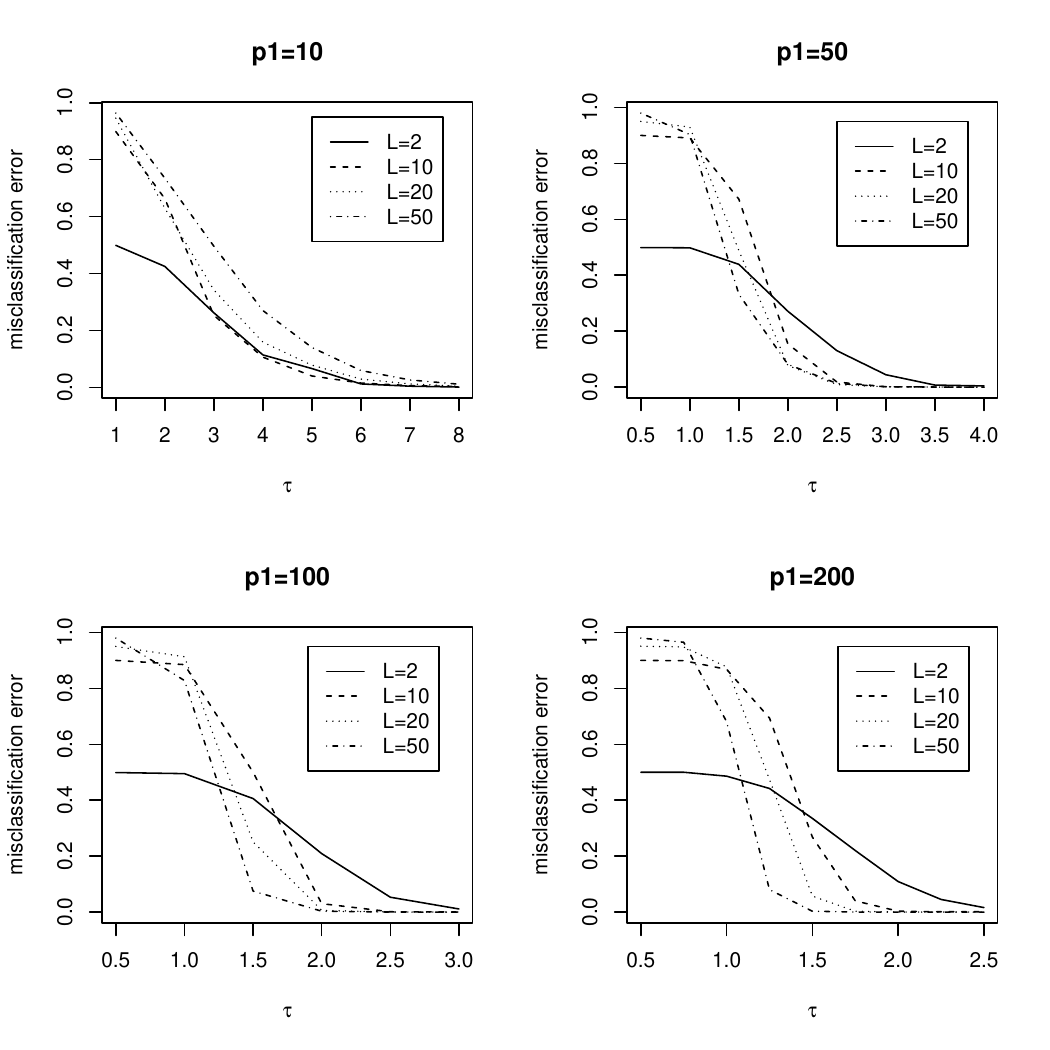}
   \caption{Average misclassification errors as  functions of $\tau$ for various combinations of $p_1$ and $L$ for Example 3.
  \label{fig:errors3} }
   \end{center}
   \end{figure}

\begin{table}
\label{tab:selection} 
\begin{center}
\begin{tabular}{|c|c||c|c|c|c||c|c|c|c|}
\hline
& & \multicolumn{4}{c||}{Example 1} & \multicolumn{4}{c|}{Example 2}\\
\hline 
$p_1$ & $\tau$ & $L=2$ & $L=10$& $L=20$& $L=50$& $L=2$ & $L=10$&$L=20$ & $L=50$\\
\hline
 10  &   1   & 1.000 &  .996 & .975  & .785  & 1.000 & 1.000 &  .978 & .788 \\
     &   2   &  .936 &  .297 & .033  & .000  &  .991 &  .592 &  .186 & .000 \\
     &   3   &  .880 &  .158 & .006  & .000  &  .898 &  .147 &  .003 & .000 \\
 50  &   1   & 1.000 &  .995 & .976  & .785  & 1.000 &  .995 &  .977 & .783 \\
     &   2   &  .975 &  .604 & .187  & .001  &  .979 &  .609 &  .172 & .001 \\
     &   3   &  .896 &  .158 & .005  & .000  &  .901 &  .146 &  .004 & .000 \\
100  &   1   & 1.000 &  .996 & .975  & .784  & 1.000 &  .996 &  .976 & .782 \\
     &   2   &  .976 &  .601 & .177  & .001  &  .981 &  .611 &  .169 & .000 \\
     &   3   &  .895 &  .149 & .005  & .000  &  .898 &  .142 &  .004 & .000 \\
200  &   1   & 1.000 &  .995 & .976  & .783  & 1.000 &  .995 &  .977 & .783 \\
     &   2   &  .975 &  .605 & .172  & .000  &  .980 &  .617 &  .175 & .000 \\
     &   3   &  .892 &  .150 & .004  & .000  &  .895 &  .150 &  .004 & .000 \\
\hline 
 \end{tabular}
\end{center}
\caption{Average proportions of false negative features 
for $p=500$ and various values of  $L$, $p_1$ and $\tau$ over $M_1 \cdot M_2=2500$ training samples.} 
\end{table}

The results of simulations indicate that for such data generating model (somewhat
different from that analyzed in the paper), the threshold $\lambda_1$ in \fr{eq:lambda1} 
(as well as $\lambda$ in \fr{eq:lambda} for the known variances) 
might be too high, especially for small values of  $\tau$. The latter led to an   over-conservative
feature selection procedure. Thus, in all simulations the feature selection procedure did not detect false positive features. 
The information on the proportions of false negative features (over the total number of significant features) 
for several combinations of $p_1$, $L$ and $\tau$ over $M_1 \cdot M_2=2500$ training samples is summarized 
in Table 1 for Example 1 and Example 2 (the results for Example 3 were similar
and we omit their presentation to save the space). In particular, Table 1 clearly shows that for small values of $\tau$ and small $L$, 
due to the over-conservative  feature selection procedure, almost not a single significant feature has been detected 
and the resulting classification  is then essentially reduced to just a pure random guess.
However, for any $\tau$ the detection rate improves as $L$ grows. The improvement rate is very fast for $\tau \geq 2$. 
Thus, for $L=50$ the vast majority of significant features were detected in spite of high level of noise.
As we have mentioned, this improves the classification precision since
weaker significant features that remained latent in coarse classification become
active and may have a strong impact with increasing $L$.

For each combination of $p_1$, $L$ and $\tau$ we 
calculated the corresponding average misclassification errors: see Figures~\ref{fig:errors1}--\ref{fig:errors3} 
for Examples~1--3, respectively. 
Figures~\ref{fig:errors1}-\ref{fig:errors3} show similar behavior for all three examples. 
For any $p_1$ and $L$ misclassification error tends to zero 
as $\tau$ increases. The decay is faster for larger $p_1$ -- the more 
significant features, the easier is classification. The figures demonstrate 
also another interesting phenomenon:  for moderate and large $p_1$, the larger   $L$, the
faster is the decay. As we have argued, this is due to the fact that the impact
of weaker significant features becomes stronger with increasing $L$.
For small $\tau$ (strong noise), misclassification errors are higher for larger number
of classes $L$. This is naturally explained by the failure of feature selection
procedure to detect significant features in this case (see comments above),
so that the resulting classification is similar to a random guess with a misclassification
error $1-1/L$ (see Figures \ref{fig:errors1}-\ref{fig:errors3}). However, as $\tau$ increases, even the 
first few detected significant features strongly improve classification precision.


\subsection{Real-data example}  \label{subsec:realdata}

We applied  feature selection techniques discussed above   to a dataset of communication signals recorded 
from South American knife fishes of the genus Gymnotus. These nocturnally active freshwater fishes generate 
pulsed electrostatic fields from electric organ discharges (EODs). The three-dimensional electrostatic 
EOD fields of Gymnotus can be summarized by two-dimensional head-to-tail waveforms recorded from 
underwater electrodes placed in front of and behind a fish. EOD waveforms vary  among species 
and are used by  genus Gymnotus in order to recognize its own kind for more productive mating and other purposes.

The data set consists of   512-dimensional vectors of  the Symmlet-4 discrete wavelet transform coefficients of signals obtained from
eight  genetically distinct species of Gymnotus ({\it G. arapaima} (G1), {\it G. coatesi} (G2),  {\it G. coropinae} (G3), 
 {\it G. curupira} (G4),  {\it G. jonasi} (G5),  {\it G. mamiraua} (G6),  {\it  G. obscurus }(G7),  {\it G. varzea} (G8)) at various stages of their development.
In particular, species were divided into six ontogenetic categories: postlarval (J0), small juvenile (J1),  large juvenile (J2),
immature adult (IA), mature male (M) and mature female (F). The EODs were recorded from 42 of 48 possible combinations 
of eight species and six  categories. There are 677 samples from 42 classes with sizes varying from 3 to 69.
The complete description of the data can be found in Crampton {\it et al.} (2011).

As it is evident from  Crampton,   Lovejoy and   Waddell (2011), there is no expectation that these 
groups should all be mutually separable: there is considerable overlaps between developmental 
stages of the same specie as well as among juveniles of different species.
For this reason, we reduced the number of classes to include only those species/categories
that might be potentially separated. In particular, we ran our feature selection and classification procedure
with the data sets comprised of 10  to 16  classes listed in the order they appear: G2-M, G4-M, G5-M, G1-F, G2-F, 
G5-F, G7-F, G8-F, G2-J1, G4-J1, G2-F, G1-J1, G7-AI, G1-F, G6-M, G7-J1.


We split  the respective data sets into training and test parts. 
For this purpose, in each class we chose at random  at most 1/3 of the total number of observations 
for validation leaving the rest of the data as training samples. Using those training samples, 
we carried out feature selection and subsequent classification of vectors in the test part of the data set.
We repeated the process 100 times for various splits and recorded the average misclassification errors and their standard errors 
for each of the cases ($L=10, 11, \dots,16$).   
Table 2 presents results of the study: the average sample sizes of train ($N_{train}$)
and test ($N_{test}$) sets for each $L$, the average number of selected significant features ($\hat{p}_1$) and average misclassification error with the corresponding standard errors.

\begin{table}
\label{tab:real_data} 
\begin{center}
\begin{tabular}{|c|c|c|c|c|}
\hline
$L$ & $N_{train}$ & $N_{test}$ & $\hat{p}_1$ & Misclassification error  \\
\hline
10 & 32 &  10  & 67.0 &  .077\; (.006) \\
11 & 38 &  13  & 68.3 &  .092\; (.006) \\ 
12 & 46 &  16  & 65.3 &  .127\; (.007) \\ 
13 & 51 &  18  & 67.6 &  .166\; (.007) \\ 
14 & 57 &  20  & 83.7 &  .149\; (.006) \\
15 & 64 &  23  & 87.4 &  .130\; (.006) \\
16 & 68 &  24  & 86.8 &  .162\; (.007) \\
\hline 
 \end{tabular}
\end{center}
\caption{The sample sizes of train ($N_{train}$)
and test ($N_{test}$) sets, the numbers of selected significant features 
($\hat{p_1}$) and misclassification errors with standard errors in brackets averaged over 100 splits for the  
Gymnotus fish data.}
 \end{table}

The table shows that when one starts with 10 well separated classes the misclassification error is initially grows 
when $L$ increases from 10 to 13. However, at $L=13$ there is a strong jump in the 
numbers of detected features and the misclassification errors again start to decrease 
when $L$  grows from 13 to 15 due to better feature selection.
For $L > 15$ the misclassification error grows again with $L$ due to poor 
separation of juvenile Gymnotus  EOD waveforms shapes.

\section{Concluding remarks} 
\label{sec:remarks}

The paper considers multi-class classification of high-dimensional normal
vectors, where the number of classes may diverge. This is 
a first attempt to rigorously  study ``large $L$, large $p$, small $n$'' classification problem. Our main goal was not to propose a novel methodology but to explore interesting phenomena arising in such a new setup.
In particular, our results indicate that 
the precision of classification can  improve as a number of classes grows. 
This is, at first glance, a somewhat counter-intuitive conclusion and has not been observed so far due to
shortage of literature on multi-class classification.
It is explained by the fact that even weaker significant
features, that might be undetected for smaller $L$, being shared across classes, can strongly contribute to
successful classification   when the   number of classes is large. 
We believe that the results of the paper motivate further investigation of ``large $L$, large $p$, small $n$'' classification 
in other, more complicated setups.

The contents of this paper can be extended in a variety of ways. 
To begin with, an extension to different covariance matrices across the classes is straightforward. 
One can also allow different supports of sparsity for different clusters and/or relax the
Gaussian assumption by considering sub-Gaussian or sub-exponential data in a similar way, 
though such generalizations will require to re-derive the corresponding conditions for correct classification.

 
\section*{Acknowledgments} 

Felix Abramovich was supported by the Israel Science Foundation (ISF), grant ISF-589/18.
Marianna Pensky   was  partially supported by National Science Foundation
(NSF), grants  DMS-1407475 and DMS-1712977. The authors would like to 
thank  Vladimir Koltchinskii and Ruth Heller for valuable remarks, and Will Crampton
for providing the data set used for the real data example. Helpful comments of the anonymous referee are gratefully acknowledged.
 

\section{Appendix}
\label{sec:appendix}

We start from recalling two lemmas of Birg\'e (2001) that will be used further in the proofs.

\begin{lemma}[\bf Lemma 8.1 of Birg\'e, 2001] \label{lem:chi2}
Let $\zeta \sim \chi^2_{k,\mu}$, $\mu>0$. Then, for any $x>0$
\be \label{eq:uptail}
P(\zeta > \mu+k + 2\sqrt{(k+2\mu)x}+2x) \leq e^{-x} 
\ee
and
\be   \label{eq:lowtail}
P(\zeta < \mu+k -2 \sqrt{(k+2\mu)x}) \leq e^{-x} 
\ee
\end{lemma}

\begin{lemma}[\bf Lemma 8.2 of Birg\'e, 2001] \label{lem:lem2bm}
Let $X$ be a random variable such that
$$
\log [E\left( e^{sX}\right)] \leq \frac{(as)^2}{1 - bs}
\quad \mbox{for} \quad 0 < s < b^{-1},
$$
where $a$ and $b$ are positive constants. Then
$$
P[X \geq 2a \sqrt{x} + bx] \leq e^{-x} \quad \mbox{for all}\ \ x > 0.
$$
\end{lemma}


\noindent
{\bf Proof of Theorem \ref{th:oracle} }\ \  
%
Note that
\be \label{eq:error1}
P(\hat{l} \neq l) = \sum_{k \neq l} P(\hat{l}=k)
\leq L_1  \max_{k \neq l}P(\hat{l}=k),
\ee
For a given $k \neq l$ define a $(2 p_1)$-dimensional random vector 
$\widetilde{\bf Y} = \begin{pmatrix} {\bf Y}^*_0-{\bf Y}^*_{l} \\ 
                  {\bf Y}^*_0-{\bf Y}^*_k
\end{pmatrix} 
$,
where the vectors ${\bf Y}^*_0, {\bf Y}^*_{l}$ and ${\bf Y}^*_k$ are defined just after (\ref{eq:trueclass}).
A straightforward calculus yields 
\be \label{eq:V}
\widetilde{\bf Y}
\sim N\left({\pmb \theta}, V\right) \quad \mbox{with} \quad
{\pmb \theta} = \begin{pmatrix} {\bf 0}_{p_1} \\ {\bf m}^*_l - {\bf m}^*_k \end{pmatrix},\quad 
V=\sigma^2 \begin{pmatrix} \rho^{-1}_l\, \Sigma^* & \Sigma^* \\
                            \Sigma^*           & \rho^{-1}_k\, \Sigma^*
            \end{pmatrix}
\ee
where $\rho_l$ is defined in (\ref{eq:rho}).  
%
%
Then, it follows from (\ref{eq:trueclass}) that
$$
P(\hat{l} = k) \leq P\left(\rho_l ({\bf Y}^*_0-{\bf Y}^*_l)^t (\Sigma^*)^{-1} ({\bf Y}^*_0-{\bf Y}^*_l) > 
\rho_k ({\bf Y}^*_0-{\bf Y}^*_k)^t (\Sigma^*)^{-1} ({\bf Y}^*_0-{\bf Y}^*_k)\right)
= P(\widetilde{\bf Y}^t A \widetilde{\bf Y} \geq 0),
$$        
where 
$$
A=\begin{pmatrix} \rho_l\, (\Sigma^*)^{-1} & 0_{p_1 \times p_1} \\
                        0_{p_1 \times p_1} &  -\rho_k\, (\Sigma^*)^{-1}
         \end{pmatrix}
$$

Consider a random variable $\xi=\widetilde{\pmb Y}^t A \widetilde{\pmb Y}$. 
Since $V^{-1}$ is a symmetric positive-definite matrix and $A$ is symmetric, they can
be simultaneously diagonalized, that is, there exists a matrix $W$, such
that $W^t V^{-1} W =I$ and $W^t A W =\Lambda$, where $\Lambda$ is a diagonal
matrix of the eigenvalues $\varphi_j,\;j=1,\ldots,2p_1$ of $R=VA$. Then, from the
known results on the distribution of quadratic forms of normal variables
(e.g., Imhof, 1961), $\xi$ can be represented as a weighted sum of independent (generally)
non-central chi-square variables, namely,
\be \label{eq:xi}
\xi =\sum_{j=1}^{2p_1}\varphi_j \chi^2_{1,\eta_j^2},
\ee
where ${\pmb \eta}$ is such that ${\pmb \theta} =W{\pmb \eta}$ with 
${\pmb \theta}$ given by (\ref{eq:V}).
By a straightforward matrix calculus, obtain 
$$
R^2=
\begin{pmatrix} 
\left(1-\rho_k \rho_l \right)~ I_{p_1}  & 0_{p_1 \times p_1}   \\ 
0_{p_1 \times p_1}  & \left(1-\rho_k \rho_l \right)~ I_{p_1} 
\end{pmatrix}
$$
and, therefore, all eigenvalues $\varphi_j,\;j=1,\ldots,2p_1,$ of   matrix $R=VA$   are 
of the forms 
\be \label{eq:eigen}
\varphi_j = \pm \varphi_*, \quad \mbox{where} \quad \varphi_*=\sqrt{1- \rho_k \rho_l},\ \ 
j=1,\ldots,2p_1
\ee


Consider now the logarithm of the moment generating function of the centered random
variable $\xi-E(\xi)$, where $\xi$ is defined in  (\ref{eq:xi}). We have 
$E\xi=\sum_{j=1}^{2p_1}\varphi_j(1+\eta_j^2)=\sum_{j=1}^{2p_1}\varphi_j \eta_j^2$,
where recall that $W{\pmb \eta} = {\pmb \theta}$.
Hence, using formula (\ref{eq:eigen}), for $s< 1/(2\varphi_*)$, we have
\be \nonumber
\begin{split}
\ln E e^{s (\xi-E\xi)} & =\sum_{j=1}^{2p_1} \frac{\eta_j^2 \varphi_j s}{1-2\varphi_j s}
-\frac{1}{2}\sum_{j=1}^{2p_1}\ln(1-2\varphi_j s) 
-s \sum_{j=1}^{2p_1}\varphi_j(1+\eta_j^2)\\
& = 
\sum_{j=1}^{2p_1} 
\left(\frac{\eta_j^2 \varphi_j s}{1-2\varphi_j s}-\eta_j^2 \varphi_j s\right)
-\frac{1}{2}\sum_{j=1}^{2p_1}\left(\ln(1-2\varphi_j s) + 2\varphi_j s\right) 
\\
& \leq \sum_{j=1}^{2p_1}\frac{2s^2 \eta_j^2 \varphi_*^2}{1-2\varphi_j s}+ 
\sum_{j=1}^{2p_1}\frac{s^2 \varphi_*^2}{1-2\varphi_j s} 
~ \leq ~ \frac{2s^2}{1-2\varphi_* s}~ \varphi_*^2 ||{\pmb \eta}||^2+
\frac{2s^2\varphi_*^2 p_1}{1-4\varphi_*^2 s^2} \\
& \leq ~ \frac{2s^2}{1-2\varphi_* s}~ \varphi_*^2 ||{\pmb \eta}||^2+
\frac{2s^2\varphi_*^2 p_1}{1-2\varphi_*  s} \\
\end{split}
\ee

Denote
$$
\Delta^2=({\bf m}^*_l-{\bf m}^*_k)^t (\Sigma^{*})^{-1} ({\bf m}^*_l-{\bf m}^*_k)
$$
Using   $W^t V^{-1} W =I$, $W^t A W =\Lambda$ and $W{\pmb \eta} = {\pmb \theta}$, one can verify that 
$ 
\varphi_*^2 ||{\pmb \eta}||^2={\pmb \eta}^t \Lambda^2 {\pmb \eta} =
{\pmb \theta}^t A V A {\pmb \theta} =  \rho_k~ \Delta^2,
$ 
where ${\pmb \theta}$ and $V$ are defined in (\ref{eq:V}).
Thus, 
$$
\ln E e^{s (\xi-E\xi)}  \leq  \frac{a^2 s^2}{1- b s},
$$
where $b=2 \varphi_*$ and
\be \nonumber
a=\sqrt{2 \rho_k\, \Delta^2 +2\varphi_*^2 p_1}
\leq \sqrt{2} \left(\sqrt{\rho_k} ~|\Delta|  +\varphi_* \sqrt{p_1}\right)
\ee
%
In addition, 
$$
E\xi={\pmb \eta}^t \Lambda {\pmb \eta}= {\pmb \theta}^t A   {\pmb \theta}
= - \rho_k~ \Delta^2
$$
A straightforward calculus shows that, under the condition (\ref{eq:as1}) of
Theorem \ref{th:oracle}, one has 
$\rho_k\, \Delta^2 \geq  2 a \sqrt{\ln (L_1/\alpha)} + b\ln (L_1/ \alpha)$.
Then, applying Lemma \ref{lem:lem2bm}, one obtains
$$
P(\xi > 0) \leq P\left(\xi \geq -\rho_k~\Delta^2+2 a \sqrt{\ln (L_1/ \alpha)} + b \ln (L_1/\alpha) \right) \leq 
\frac{\alpha}{L_1}
$$
that, together with (\ref{eq:error1}), complete the proof.


\medskip

 
\noindent
{\bf Proof of Theorem  \ref{th:pi1} }\ \  
Let $\hat{p}_{01}=\sum_{j=1}^p I\{\hat{x}_j=1 \mid x_j=0\}$ and $\hat{p}_{11}=\sum_{j=1}^p I\{\hat{x}_j=1 \mid x_j=1\}$ 
be the numbers of erroneously and truly identified significant features respectively, where obviously $\hat{p}_{01}$ and $\hat{p}_{11}$ are
independent, and $\hat{p}_{01}+\hat{p}_{11}=\hat{p}_1$. 
Note that
$$
P(\hat{x} \neq x) \leq P(\hat{p}_{01} > 0) + P(\hat{p}_{11} < p_1)
$$

Recall that for $x_j=0$, the corresponding $\zeta_j \sim \chi^2_{L_1}$. 
Let $u_j,\;j=1,\ldots,p_0$ be any, possibly correlated,  $\chi^2_{L_1}$ random variables. Then,
\be \nonumber
P(\hat{p}_{01}> 0)=P\left(\max_{1 \leq j \leq p_0} u_j > \lambda\right) 
\leq p~ P\left(u_j>L_1+2\sqrt{L_1 \ln (2p/\alpha)}+2 \ln(2p/\alpha) \right)
\ee
Apply  Lemma \ref{lem:chi2} for the particular case
$\mu=0$ to obtain
$$
P\left(u_j > L_1+2\sqrt{L_1 \ln(2p/\alpha)}+2 \ln(2p/\alpha)\right) \leq
\frac{\alpha}{2p},
$$
so that $P(\hat{p}_{01}> 0) \leq \alpha/2$.
Similarly, let $\mu_*=\min_{1 \leq j \leq p_1} \mu_j=
\min_{1 \leq j \leq p_1} \sigma^{-2}_j~\sum_{l=1}^L n^{(1)}_l\beta_{lj}^2$ and 
consider any, possibly  correlated, non-central chi-squared variables
$v_j \sim \chi^2_{L_1;\mu_*},\;j=1,\ldots, p_1$. We have
\be \nonumber
P(\hat{p}_{11} < p_1) \leq P\left(\min_{1 \leq j \leq p_1}
v_j \leq  \lambda\right)  \leq p~ P\left(v_j < \lambda\right)
\ee
A straightforward calculus shows that, under the condition (\ref{eq:as2}) on $\mu_*$, 
one has $\mu_* +L_1-2\sqrt{(L_1+2\mu_*)\ln (2p/\alpha)} > \lambda$. Thus, Lemma \ref{lem:chi2}
yields $P(v_j < \lambda) \leq \alpha/(2p)$ and, therefore,
$P(\hat{p}_{11} < p_1) \leq \alpha/2$, which completes the proof.

\medskip

 
\noindent
{\bf Proof of Theorem  \ref{th:pi2} }\ \  
%
We start with the following lemma:
\begin{lemma} \label{lem:sigmahat}
$$
P\left(\max_{1 \leq j \leq p} \left|\hat{\sigma}_j^2/\sigma^2 -1\right| \leq \kappa \right) \geq 1-\alpha,
$$
where $\kappa$ was defined in (\ref{eq:kappa}).
\end{lemma}

Let $\cA$ be the event $\{\max_{1 \leq j \leq p} \left|\hat{\sigma}_j^2/\sigma^2 -1\right| \leq \kappa\}$ and 
$I_{\cA}$ its indicator. By Lemma \ref{lem:sigmahat}, 
\be \label{eq:A1}
P(\hat{x} \neq x) \leq P\left((\hat{x} \neq x) I_{\cA}\right) + \alpha,
\ee
where
\be \label{eq:A2}
P\left((\hat{x} \neq x)I_{\cA}\right) \leq P\left((\hat{p}_{01}>0)I_{\cA}\right)+P\left((\hat{p}_{11}<p_1)I_{\cA}\right)
\ee
Let $\hat{\zeta}_j=\hat{\sigma}^{-2}_j~\sum_{l=1}^L n_l^{(1)} (\bar{Y}_{lj}^{(1)} -\bar{Y}_{\cdot j}^{(1)})^2$. Then, on the event $\cA$
$$
P\left((\hat{\zeta}_j > \lambda_1)I_{\cA} \mid x_j=0\right)=P\left(\left(u_j > \lambda_1\,  \hat{\sigma}^2_j/\sigma^2_j \right)
I_{\cA}\right) \leq P(u_j > \lambda)
$$
where $u_j \sim \chi^2_{L_1},\;j=1,\ldots,p_0$. Hence, following the arguments of Theorem \ref{th:pi1}, by Lemma \ref{lem:chi2} 
\be \label{eq:A3}
P\left((\hat{p}_{01}>0)I_{\cA}\right) \leq P\left((\max_{1 \leq j \leq p} \hat{\zeta}_j > \lambda_1)I_{\cA} \mid x_j=0\right) \leq P(\max_{1 \leq j \leq p_0} u_j > \lambda) \leq \frac{\alpha}{2}
\ee
Similarly,
$
P\left((\hat{\zeta}_j < \lambda_1)I_{\cA} \mid x_j=1\right) \leq P\left(v_j < \lambda_1 (1+\kappa)\right)
$
where $v_j \sim \chi^2_{L_1;\mu_*}, j=1,\ldots,p_1$. Then, under the condition (\ref{eq:as2})
of the theorem, Lemma \ref{lem:chi2} yields
\be \label{eq:A4}
P\left((\hat{p}_{11} < p_1)I_{\cA}\right) \leq P\left(\min_{1 \leq j \leq p_1} v_j \leq \lambda_1
(1+\kappa)\right) \leq \frac{\alpha}{2}  
\ee
Combination of (\ref{eq:A1})-(\ref{eq:A4}) completes the proof.


\medskip

 
\noindent
{\bf Proof of Theorem  \ref{th:main1} }\ \  
%
Assume that  ${\bf Y}_0$ is from the $l$-th class. From (\ref{eq:error}) we have
$ 
P(\hat{l} \neq l) \leq P(\hat{l}\neq l\mid\hat{x}=x) + P(\hat{x} \neq x),
$ 
where $P(\hat{x} \neq x) \leq 2\alpha$ by Theorem \ref{th:pi2}.
Consider a set $\Omega = \{ \omega:\  \hat{x}=x \}$ with $P(\Omega) \geq 1- \alpha$.
In order to bound above  $P(\hat{l} \neq l \mid \hat{x}=x)$ we assume that $\omega \in \Omega$. 
We will use the following two lemmas:


\begin{lemma} \label{lem:sigmainv}
If $||\widehat{\Sigma^*}-\Sigma^*|| \leq  \lambda_{\min}(\Sigma^*)/2$, then
$
||(\widehat{\Sigma^*})^{-1}-(\Sigma^*)^{-1}|| \leq  2~ \lambda^{-2}_{\min}(\Sigma^*)\, ||\widehat{\Sigma^*}-\Sigma^*||
$
\end{lemma}


\begin{lemma} \label{lem:sigma}
Under the condition (\ref{eq:p1N}),\ \   
$
P\left(||\widehat{\Sigma^*}-\Sigma^*|| \leq  \lambda_{\max}(\Sigma^*)\sqrt{\frac{C_1 p_1}{N_2}} \right) \geq 1-2\alpha
$
\end{lemma}


From Lemma \ref{lem:sigmainv} and Lemma \ref{lem:sigma} it follows  that
under (\ref{eq:p1N}), 
\be \label{eq:sigmainv}
P\left(||(\widehat{\Sigma^*})^{-1}-(\Sigma^*)^{-1}|| \leq \gamma_{p_1,N_2}
\right) \geq 1-2\alpha
\ee
where $\gamma_{p_1,N_2}$ is defined in \eqref{eq:betap1n}.
Furthermore, for any $1 \leq k \leq L$, 
\be \label{eq:a4}
\frac{({\bf Y}^*_0-\bar{\bf Y}^*_k)^t \left((\widehat{\Sigma^*})^{-1}-(\Sigma^*)^{-1}\right)({\bf Y}^*_0-\bar{\bf Y}^*_k)}
{({\bf Y}^*_0-\bar{\bf Y}^*_k)^t (\Sigma^*)^{-1} ({\bf Y}^*_0-\bar{\bf Y}^*_k)}
\leq ||\Sigma^*\left((\widehat{\Sigma^*})^{-1}-(\Sigma^*)^{-1}\right)|| \leq \tau_2 ||(\widehat{\Sigma^*})^{-1}-(\Sigma^*)^{-1})||
\ee
Since the sample mean and the sample covariance matrix are independent in the case of the normal distribution, 
inequalities  (\ref{eq:sigmainv}) and (\ref{eq:a4}) imply that with probability at least $1-2\alpha$
\be \nonumber
\begin{split} 
 & \rho_l \, ({\bf Y}^*_0-\bar{\bf Y}^{(2)*}_l)^t (\widehat{\Sigma^*})^{-1}({\bf Y}^*_0-\bar{\bf Y}^{(2)*}_l)- 
\rho_k \, ({\bf Y}^*_0-\bar{\bf Y}^{(2)*}_k)^t (\widehat{\Sigma^*})^{-1}({\bf Y}^*_0-\bar{\bf Y}^{(2)*}_k) \\
&= \rho_l \, ({\bf Y}^*_0-\bar{\bf Y}^{(2)*}_l)^t (\Sigma^*)^{-1}({\bf Y}^*_0-\bar{\bf Y}^{(2)*}_l)- 
\rho_k \, ({\bf Y}^*_0-\bar{\bf Y}^{(2)*}_k)^t (\Sigma^*)^{-1}({\bf Y}^*_0-\bar{\bf Y}^{(2)*}_k) \\
& +
\rho_l \, ({\bf Y}^*_0-\bar{\bf Y}^{(2)*}_l)^t \left((\widehat{\Sigma^*})^{-1}-(\Sigma^*)^{-1}\right)({\bf Y}^*_0-\bar{\bf Y}^{(2)*}_l)- 
\rho_k \, ({\bf Y}^*_0-\bar{\bf Y}^{(2)*}_k)^t \left((\widehat{\Sigma^*})^{-1}-(\Sigma^*)^{-1}\right)({\bf Y}^*_0-\bar{\bf Y}^{(2)*}_k) \\
& \leq \rho_l(1+\gamma_{p_1,N_2}) \, ({\bf Y}^*_0-\bar{\bf Y}^{(2)*}_l)^t (\Sigma^*)^{-1}({\bf Y}^*_0-\bar{\bf Y}^{(2)*}_l)-\rho_k(1-\gamma_{p_1,N_2}) \, 
({\bf Y}^*_0-\bar{\bf Y}^{(2)*}_k)^t (\Sigma^*)^{-1}({\bf Y}^*_0-\bar{\bf Y}^{(2)*}_k)
\end{split}
\ee
Define $\rho'_l=\rho_l(1+\gamma_{p_1,N_2})$ and $\rho'_k=\rho_k(1-\gamma_{p_1,N_2})$. In particular,
note that $\rho'_l \rho'_k = \rho_l \rho_k (1-\gamma_{p_1,N_2}^2)$. 
Repeating the proof of Theorem \ref{th:oracle}  but with
$\rho'_l$ and $\rho'_k$ and under the stronger condition
(\ref{eq:as11}), obtain $P(\hat{l} \neq l \mid \hat{x}=x) \leq 2\alpha$ that, 
together with (\ref{eq:error}) and $P(\hat{x} \neq x) \leq 2\alpha$, completes the proof.

\medskip

 
\noindent
{\bf Proof of Lemma \ref{lem:sigmahat} }\ \  
Note that $\sigma_j^{-2}(N_1-L)\hat{\sigma}_j^2  \sim \chi^2_{N_1-L}$ and apply
Lemma \ref{lem:chi2} to obtain 
$P(|\hat{\sigma}_j^2/\sigma^2 -1| \geq \kappa) \leq  \alpha/p$ 
for all $j=1,\ldots,p$ and, therefore,
$P \lkr \max_{1 \leq j \leq p} |\hat{\sigma}_j^2/\sigma^2 -1| \geq \kappa \rkr \leq  \alpha$ 

\ignore{ 
$$
P\left(\left|\frac{\hat{\sigma}_j^2}{\sigma^2}-1\right| \geq \kappa \right) \leq \frac{\alpha}{p}
$$
for all $j=1,\ldots,p$ and, therefore,
$$
P\left(\max_{1 \leq j \leq p} \left|\frac{\hat{\sigma}_j^2}{\sigma^2}-1\right| \geq \kappa \right) \leq \alpha
$$
}

\medskip

 
\noindent
{\bf Proof of Lemma \ref{lem:sigmainv} }\ \  
%
Under the condition of the lemma we have
$$
||(\widehat{\Sigma^*})^{-1}||^{-1}=\min_{||\bx||=1} \bx^t \widehat{\Sigma^*} \bx \geq \min_{||\bx||=1} \bx^t \Sigma^* \bx - 
\max_{||\bx||=1} \bx^t(\widehat{\Sigma^*}-\Sigma^*)\bx \geq  \lambda_{\min}(\Sigma^*)/2
$$
and, therefore,
$$
||(\widehat{\Sigma^*})^{-1}-(\Sigma^*)^{-1}|| \leq ||(\widehat{\Sigma^*})^{-1}|| \cdot ||\widehat{\Sigma^*}-\Sigma^*|| 
\cdot ||(\Sigma^*)^{-1}|| \leq  2 \lambda^{-2}_{\min}(\Sigma^*)\, ||\widehat{\Sigma^*}-\Sigma^*|| 
$$

\medskip

 
\noindent
{\bf Proof of Lemma \ref{lem:sigma} }\ \  
%
Define ${\bf Z}_{il}= \lkr{\bf Y}^*_{il}\rkr^{(2)} -{\bf m}_l^* \sim N({\bf 0}_{p_1},\Sigma^*),\;i=1,\ldots,n_l^{(2)};\;l=1,\ldots,L$. 
The sample covariance matrix is translation invariant and, therefore, 
$$
\widehat{\Sigma^*}=\frac{1}{N_2}\sum_{l=1}^L \sum_{i=1}^{n_l^{(2)}} ({\bf Z}_{il}-\bar{\bf Z}_l)({\bf Z}_{il}-\bar{\bf Z}_l)^t=
\frac{1}{N_2} \sum_{l=1}^L \sum_{i=1}^{n_l^{(2)}} {\bf Z}_{il}{\bf Z}_{il}^t-
\frac{1}{N_2} \sum_{l=1}^L n_l^{(2)}~ \bar{\bf Z}_l \bar{\bf Z}_l^t= S_1 - S_2
$$
Thus,
\be \label{eq:m1}
||\widehat{\Sigma^*}-\Sigma^*|| \leq ||S_1-\Sigma^*||+||S_2||
\ee
By Remark 5.51 of Vershynin (2012), under the conditions of the lemma there
exists an absolute constant $C_0$ such that
\be \label{eq:C1}
P\left(||S_1-\Sigma^*|| \leq \tau_2 \sqrt{\frac{C_0 p_1}{N_2}}\right) \geq 1-\alpha
\ee
Consider now $S_2$. Define the $p_1 \times L$-dimensional matrix  $\bar{Z}$ with columns 
$\bar{\bf Z}_l$, $l=1, \cdots, L$ and the  diagonal matrix 
$D = \mbox{diag}(\sqrt{n_1^{(2)}}, \cdots, \sqrt{n_L^{(2)}})$. It is easy to see that
$S_2 = N^{-1}\, (\bar{Z} D) (\bar{Z} D)^t$ and that matrix $\Xi = (\Sigma^*)^{-1/2} \bar{Z} D$
has i.i.d. $N(0,1)$ entries. Indeed, columns 
${\bf \Xi}_l = \sqrt{n_l^{(2)}} (\Sigma^*)^{-1/2}\, \bar{\bf Z}_l$ of matrix $\Xi$
are independent with $\Cov ({\bf \Xi}_l) = I_{{p_1}}$. Hence,
$$
\| S_2\| = N_2^{-1}\, \|\bar{Z} D \|^2 = N_2^{-1}\, \| \sqrt{\Sigma^*}\, \Xi\|^2 \leq 
 N_2^{-1}\, \lambda_{\max}(\Sigma^*) \| \Xi\|^2.
$$
Then, by Corollary 5.35 of Vershynin (2012) 
$$
P\left(||S_2|| \leq   N_2^{-1}\,\lambda_{\max}(\Sigma^*)  \left(\sqrt{p_1}+\sqrt{L}+\sqrt{2\ln(2/\alpha)}\right)^2 \right)
\geq 1-\alpha
$$
that, under (\ref{eq:p1N}), yields 
\be \label{eq:C2}
P\left(||S_2|| \leq 9  \lambda_{\max}(\Sigma^*) N_2^{-1}\, p_1 \right) \geq 1-\alpha
\ee
Combination of  (\ref{eq:m1})-(\ref{eq:C2}) completes the
proof with $C_1=\max(\sqrt{C_0},9)$.



\end{document}